\catcode`\@=11

\magnification=1200
\baselineskip=14pt

\pretolerance=500    \tolerance=1000 \brokenpenalty=5000

\catcode`\;=\active
\def;{\relax\ifhmode\ifdim\lastskip>\z@
\unskip\fi\kern.2em\fi\string;}

\overfullrule=0mm

\catcode`\!=\active
\def!{\relax\ifhmode\ifdim\lastskip>\z@
\unskip\fi\kern.2em\fi\string!}

\catcode`\?=\active
\def?{\relax\ifhmode\ifdim\lastskip>\z@
\unskip\fi\kern.2em\fi\string?}

\frenchspacing

\newif\ifpagetitre            \pagetitretrue
\newtoks\hautpagetitre        \hautpagetitre={ }
\newtoks\baspagetitre         \baspagetitre={1}

\newtoks\auteurcourant        \auteurcourant={M. Laurent   }
\newtoks\titrecourant
\titrecourant={ On transfer  inequalities in Diophantine Approximation }

\newtoks\hautpagegauche       \newtoks\hautpagedroite
\hautpagegauche={\hfill\sevenrm\the\auteurcourant\hfill}
\hautpagedroite={\hfill\sevenrm\the\titrecourant\hfill}

\newtoks\baspagegauche       \baspagegauche={\hfill\rm\folio\hfill}

\newtoks\baspagedroite       \baspagedroite={\hfill\rm\folio\hfill}

\headline={
\ifpagetitre\the\hautpagetitre
\global\pagetitrefalse
\else\ifodd\pageno\the\hautpagedroite
\else\the\hautpagegauche\fi\fi}

\footline={\ifpagetitre\the\baspagetitre
\global\pagetitrefalse
\else\ifodd\pageno\the\baspagedroite
\else\the\baspagegauche\fi\fi}

\def\date{\ {\the\day}\
\ifcase\month\or Janvier\or F\'evrier\or Mars\or Avril
\or Mai \or Juin\or Juillet\or Ao\^ut\or Septembre
\or Octobre\or Novembre\or D\'ecembre\fi\
{\the\year}}

\def\up#1{\raise 1ex\hbox{\sevenrm#1}}

\def\cqfd{\unskip\kern 6pt\penalty 500
\raise -2pt\hbox{\vrule\vbox to 10pt{\hrule width 4pt
\vfill\hrule}\vrule}\par\medskip}

\def\section#1{\vskip 7mm plus 20mm minus 1.5mm\penalty-50
\vskip 0mm plus -20mm minus 1.5mm\penalty-50
{\bf\noindent#1}\nobreak\smallskip}

\def\subsection#1{\medskip{\bf#1}\nobreak\smallskip}

\def\displaylinesno #1{\dspl@y\halign{
\hbox to\displaywidth{$\@lign\hfil\displaystyle##\hfil$}&
\llap{$##$}\crcr#1\crcr}}

\def\ldisplaylinesno #1{\dspl@y\halign{
\hbox to\displaywidth{$\@lign\hfil\displaystyle##\hfil$}&
\kern-\displaywidth\rlap{$##$}
\tabskip\displaywidth\crcr#1\crcr}}

\def\hfl#1#2{\smash{\mathop{\hbox to 12 mm{\rightarrowfill}}
\limits^{\scriptstyle#1}_{\scriptstyle#2}}}

%%%%%%%%%%%%%%%%%%%%%%%%%%%%%%%%%%%%%%%%%%%%%%%%%%%%%%%%%%%%%%%%%
\catcode`\@=12

\def\bP{{\bf P}}
\def\bQ{{\bf Q}}
\def\bx{{\bf  x}}
\def\bX{{\bf X}}
\def\by{{\bf y}}

\def\bz{{\bf z}}

\def\bR{{\bf R}}
\def\bZ{{\bf Z}}

\def\th{{\theta}}
\def\Th{{\Theta}}

\def\bZ{{\bf Z}}

\def\proof{\bigskip\noindent{\it Proof.}\ }

\def\and{\quad\hbox{and}\quad}

\def\om{{\omega}}
\def\omc{{\hat{\omega}}}

%\dimen1=\ht1 \advance\dimen1 by 2pt \dimen2=\dp1 \advance\dimen2 by 2pt
%\setbox1=\hbox{\vrule height\dimen1 depth\dimen2\box1\vrule}%
%\setbox1=\vbox{\hrule\box1}%
%\advance\dimen1 by .4pt \ht1=\dimen1
%\advance\dimen2 by .4pt \dp1=\dimen2 \box1\relax

\def\M{\mathop{\rm M\kern 1pt}\nolimits}
\def\h{\mathop{\rm h\kern 1pt}\nolimits}

\def\romain#1{\uppercase\expandafter{\romannumeral #1}}

\def\card{\mathop{\rm Card\kern 1.3 pt}\nolimits}
\def\deg{\mathop{\rm deg\kern 1pt}\nolimits}
\def\det{\mathop{\rm det\kern 1pt}\nolimits}

\def\h{\mathop{\rm h\kern 1pt}\nolimits} \long\def\forget#1\endforget{}

\def\og{\leavevmode\raise.3ex\hbox{$\scriptscriptstyle 
\langle\!\langle\,$}}
\def\fg{\leavevmode\raise.3ex\hbox{$\scriptscriptstyle
\!\rangle\!\rangle\,\,$}}

%%%%%%%%%%%%%%%%%%%%%%%%%%%%%%%%%%%%%%%%%%%%%%%%%%%%%
%%%%%%%%%%%% Biblio %%%%%%%%%%%%%%%¢

\def\BuLaB{1}
\def\BuLaC{2}
\def\Bur{3}
\def\Cas{4}
\def\DaSc{5}
\def\Dys{6}
\def\Gro{7}
\def\JarA{8}
\def\JarC{9}
\def\JarD{10}
\def\KhiA{11}
\def\KhiB{12}
\def\LauA{13}
\def\LauB{14}
\def\Mah{15}
\def\Phi{16}
\def\Rot{17}
\def\RoyA{18}
\def\RoyC{19}
\def\RoyB{20}
\def\SchA{21}
\def\SchB{22}

\centerline{}

\vskip 4mm

\centerline{
\bf  On transfer inequalities in Diophantine Approximation}

\vskip 8mm
\centerline{ by Michel L{\sevenrm AURENT}\footnote{}{\rm
} }

\vskip 9mm

\noindent {\bf Abstract --}  Let $\Th$ be a point in $\bR^n$. We split the classical Khintchine's Transference Principle
into $n-1$ intermediate estimates which connect exponents $\om_d(\Th)$ measuring the sharpness
 of the  approximation to $\Th$ by linear rational varieties
of dimension $d$,  for  $0\le d \le n-1$.  We also review old and recent results  related to  these $n$ exponents.

\vskip15mm

\section
{1. Introduction and results.}
We revisit in this note some results of homogeneous transfer in Diophantine Approximation. Let us first recall the classical
{\it Khintchine's Transference Principle} in its primary  form [\Cas , \KhiB]. Let $n$ be a positive  integer 
and let $\Th =(\th_1, \dots , \th_n)$ be a point in $\bR^n$.
We shall assume in all the forthcoming statements that the real numbers $1, \th_1, \dots , \th_n$ are linearly independent over 
the field $\bQ$ of rational numbers. Khintchine's Transference Principle relates the sharpness of the rational simultaneous 
approximation to $\th_1, \dots , \th_n$
with the measure of linear independence over $\bQ$ of $1,\th_1, \dots , \th_n$. Let us first quantify these notions.

\proclaim 
Definition 1. We denote respectively by $\om_0(\Th)$ and $\om_{n-1}(\Th)$ (the meaning  of the indices $0$ and $n-1$ will be explained
 afterwards) the supremum, possibly infinite,  of the real numbers $\om$ for which there exist infinitely many 
integer $(n+1)$-tuples $(x_0, \dots  , x_n)$ satisfying respectively the inequation
$$
\max_{1\le i\le n} |x_0 \th_i - x_i |  \le \Big( \max_{0\le i\le n} |  x_i | \Big)^{-\om}
\quad {\sl or} \quad 
 |x_0  + x_1 \th_1 + \cdots +  x_n \th_n |  \le \Big( \max_{0\le i\le n} |  x_i  | \Big)^{-\om}.
$$ 

Now we can state Khintchine's Transference Principle  as follows:

\proclaim 
Theorem 1. The inequalities 
$$
{\om_{n-1}(\Th)\over (n-1)\om_{n-1}(\Th) +n} \le \om_0(\Th) \le { \om_{n-1}(\Th) -n+1\over n}
\leqno{(1)}
$$
hold for any point $\Th$ in $\bR^n$ with $1, \th_1, \dots , \th_n$ linearly independent over $\bQ$. 
Moreover,  both  inequalities are optimal.

It is convenient to view $\bR^n$ as a subset of $\bP^n(\bR)$
via the usual embedding $(x_1,\dots ,x_n)\mapsto (1,x_1,\dots ,x_n)$. In the sequel, we shall  identify $\Th$
 with its image in $\bP^n(\bR)$. Let us introduce for each integer $d$ with 
 $0 \le d \le n-1$, an exponent $\om_d(\Th)$ which measures  the approximation to the point $\Th \in \bP^n(\bR)$
 by rational linear projective subvarieties of dimension $d$, in term of their height.  Denote by d the projective distance 
 on $\bP^n(\bR)$ (it will be defined in Section 2 below; notice  however that the use of any locally 
 equivalent distance, such as the distance associated
 to the supremum norm on $\bR^n$, would lead to the same exponents).
 For any real linear subvariety $L$ of $\bP^n(\bR)$, we denote by 
 $$
 {\rm d} ( \Th , L) = \min_{ P\in L} {\rm d} (\Th , P)
 $$
 the minimal distance between $\Th$ and the real  points $P$ of $L$. 
 When $L$ is rational over $\bQ$, we indicate  moreover by $H(L)$ its height, that is the Weil height of any system of Pl\" ucker coordinates 
 of $L$. The Weil height is normalized by using the Euclidean norm at the Archimedean place of $\bQ$.

 \proclaim
Definition 2. Let $d$ be an integer with $0\le d \le n-1$.  We denote by $\om_d(\Th)$   the supremum 
 of the real numbers $\om$ for which there exist infinitely many rational linear subvarieties $L \subset \bP^n(\bR)$
 such that 
$$
\dim (L) = d \and {\rm d}(\Th , L) \le H(L)^{-1-\om}.
$$ 

Definitions 1 and 2 are consistent, since ${\rm d}(\Th, L)$ compares respectively with 
$$
\max_{1\le i\le n} \Big| \th_i -{x_i\over x_0}\Big|  \and {| y_0 + y_1 \th_1 + \cdots + y_n \th_n| \over {\displaystyle\max_{0 \le i \le n} | y_i|}}
$$
when $L$ is either the rational point (case $d=0$) with homogeneous  coordinates $(1,x_1/x_0, \dots , x_n/x_0)$,  or the hyperplane 
(when $d=n-1$) with homogeneous equation $y_0 X_0 + \cdots + y_n X_n=0$. 

W. Schmidt was the first to investigate in [\SchA] the properties of these exponents $\om_d(\Th)$.
In fact, he did  not introduce them explicitely; however  Theorems 9--18 of [\SchA] provide us with various  relations
 between more general families of exponents,  which coincide with the exponents $\om_d(\Th)$ in our more limited framework. 
 Note that in the setting of  [\SchA], 
 the point $\Th$, which belongs here to $\bP^n(\bR) = \bP(\bR^{n+1})$,  has to be  replaced by  the associated line in $\bR^{n+1}$. 
 The  consideration of  real linear subspaces $\Th$ of $\bR^{n+1}$
 with arbitrary dimension, gives then rise to various angles and consequently to  further exponents.  
 See Section  4 of [\BuLaC] for a projective formulation similar to the present one. Here is our main result
 which  slightly improves on earlier inequalities due to W. Schmidt [\SchA]. 
 
 \proclaim
 Theorem 2. For any point $\Th$ with $1,\th_1, \dots , \th_n$ linearly independent over $\bQ$, 
 and for any integer $d$ with $1 \le d \le n-1$, we have the estimate
 $$
 { d \, \om_d(\Th)  \over \om_d(\Th) +  d + 1} \le
 \om_{d-1}(\Th) \le
 { (n-d) \om_d(\Th) - 1 \over n-d+1}.
 $$
 
 Following the terminology of [\SchA], we split the above estimate into two parts; namely   the {\it Going-up transfer}
 $$
 \om_{d+1}(\Th) \ge { (n-d) \om_d(\Th) +1 \over n-d-1}, \quad 0 \le d \le n-2,
 \leqno{(2)}
 $$
 and the {\it Going-down transfer}
 $$
 \om_{d-1}(\Th) \ge { d \, \om_d(\Th)  \over \om_d(\Th) +  d + 1}, \quad 1\le d \le n-1.
 \leqno{(3)}
 $$
The Going-up transfer  (2) refines Theorem 11 of [\SchA] in our restricted situation. It will be established  in Section 3. 
The Going-down transfer (3) is a very special case of Theorem 10  from  [\SchA].
 We  shall not prove it  here again, although a similar argumentation
using the inner product instead of the wedge product, could be employed. Note that  exterior and inner products
are algebraic translations of the geometrical operations of join and intersection for  subvarieties in a projective space.

An interesting feature of Theorem 2 is that it covers Khintchine's Transference Principle. Iterating  respectively (2) and (3), we easily find
by induction on the integer $k$ the  lower bounds
$$
 \om_{d+k}(\Th) \ge { (n-d) \om_d(\Th) +k \over n-d-k}, \quad  0 \le d \le n-2, \, 1 \le k \le n-d-1,
 \leqno{(4)}
 $$
 and
 $$
  \om_{d-k}(\Th) \ge { (d-k+1) \om_d(\Th)  \over k \, \om_d(\Th) +  d + 1}, \quad 1 \le d \le n-1,  \, 1 \le k \le d.
  \leqno{(5)}
  $$
Thus,  selecting  respectively in (4) and (5)   the pairs
$(d,k)= (0,n-1)$  and $(d,k)= (n-1,n-1)$, we obtain (1).

Jarn\'\i k [\JarA, \JarC] has established that Khintchine's Transference Principle is best possible. Since it is a formal corollary of all the inequalities
occurring in (2) and (3), it follows that any of these inequalities is optimal. 
Let us call {\it spectrum} of the function $\om_d$, the set of values taken by the exponents $\om_d(\Th)$, when $\Th$ ranges over 
$\bR^n$, with $1, \th_1, \dots , \th_n$ linearly independent over $\bQ$. Then, Theorem 2 implies the following 

\proclaim
Corollary. For any integer $d$ with $0 \le d \le n-1$, we have the lower bound
$$
\om_d(\Th)  \ge { d+1\over n-d}, \leqno{(6)}
$$
and equality holds in (6) for almost all $\Th$ with respect to Lebesgue measure on $\bR^n$. The spectrum
of the function $\om_d$ is equal to the whole interval $[(d+1)/(n-d), +\infty]$.

\proof We deduce from (4) and (5) the lower bounds
$$
 \om_{d}(\Th) \ge { n \, \om_0(\Th) +d \over n-d}
 \and
  \om_{d}(\Th) \ge { (d+1) \om_{n-1}(\Th)  \over (n-d-1)\om_{n-1}(\Th) +  n}.
  $$
 By Dirichlet's Box Principle, we know that $\om_0(\Th) \ge 1/n$ and $\om_{n-1}(\Th) \ge n$.
 Then, each of the two preceding inequalities imply (6). Moreover, (5)  and (4) also provide  us with
 the converse estimates
 $$
  \om_{0}(\Th) \ge {  \om_d(\Th)  \over d \, \om_d(\Th) + d +1}
 \and
  \om_{n-1}(\Th) \ge { (n-d) \om_{d}(\Th)  +n-d-1},
  $$
so that $ \om_{0}(\Th)> 1/n$ and $\om_{n-1}(\Th) > n$ whenever $\om_d(\Th)  > ( d+1)/( n-d)$. 
Now, the easy part of the classical Khintchine-Groshev Theorem [\Gro, \KhiA] shows that the set of points $\Th$ for which
$\om_{0}(\Th) > 1/n$, or equivalently $ \om_{n-1}(\Th) > n$, has Lebesgue measure 0 in $\bR^n$. 
Thus, equality holds in (6) for almost all $\Th$.

As for the spectrum of the functions $\om_d$, notice that  all  lower bounds (4)  turn out to be  equalities 
$$
 \om_{d+k}(\Th) =  { (n-d) \om_d(\Th) +k \over n-d-k} , \quad  0\le d \le n-2, 1 \le k \le n-d-1,
  $$
 whenever 
 $
 \om_{n-1}(\Th) 
 $
 has  minimal value $n \, \om_0(\Th)  +n -1$, with regard to Khintchine's Transference Principle (1). 
 Now,  Jarn\'\i k [\JarC] has established  that for any $w$ with $1/n \le w \le + \infty$, there exists
 a point $\Th\in \bR^n$ such that
 $$
 \om_0(\Th) = w \and \om_{n-1}(\Th) = n w +n-1.
 $$
 For such a point $\Th$, we then have
 $$
 \om_{d}(\Th) ={ n w +d \over n-d}, \quad 0 \le d \le n-1.
 $$
 It follows that, for any $0\le d \le n-1$, the spectrum of the function $\om_d$ coincides with the interval $[(d+1)/(n-d), +\infty]$.
\cqfd

\noindent
{\bf Remark. } Similarly, the inequalities (5) become equalities when
$$
\om_0(\Th) = { \om_{n-1}(\Th) \over (n-1) \om_{n-1}(\Th) +n}.
$$
 Using now the construction of [\JarA],  we obtain  in this way a second family
$$
\om_d(\Th) = { (d+1) w \over (n-d-1)w +n}, \quad 0 \le d \le n-1,
$$
of $n$-tuples $(\om_0(\Th), \dots , \om_{n-1}(\Th))$ indexed by the real parameter $w\ge n$. We are thus  lead to address
the following
\proclaim
Problem. Find the spectrum in $(\bR \cup \{+ \infty\})^n$ of the $n$-tuples
$$
\Big(\om_0(\Th), \dots , \om_{n-1}(\Th)\Big), 
$$
when $\Th$ ranges over $\bR^n$,  with $1, \th_1, \dots , \th_n$ linearly independent over $\bQ$.

Necessary assumptions are provided by the inequalities (4) and (5), together with the lower bounds (6).
Are these conditions sufficient ? It holds true for $n=2$, as follows from  our work [\LauB] whose content is described
below.

\bigskip

Let us now display   a different kind of refinement of Khintchine's Transference Principle, that has  recently been 
pointed out by Y. Bugeaud and M. Laurent in  Theorem 8 of [\BuLaC]. Following the  ``hat'' notations of [\BuLaB],  introduce first {\it uniform}
analogues of the above-mentioned exponents $\om_0(\Th)$ and $\om_{n-1}(\Th)$.

\proclaim 
Definition 3. We denote respectively by $\omc_0(\Th)$ and $\omc_{n-1}(\Th)$ the supremum of the real numbers $\om$ 
such that for all sufficiently large real number $X$,  there exists a non-zero
integer $(n+1)$-tuple $(x_0, \dots  , x_n)$,  with supremum norm
$
\displaystyle\max_{0\le i \le n} | x_i | \le X
$, which  satisfies respectively the inequation
$$
\max_{1\le i\le n} |x_0 \th_i - x_i |  \le X^{-\om}
\quad {\sl or} \quad 
 |x_0  + x_1 \th_1 + \cdots +  x_n \th_n |  \le X^{-\om}.
$$ 

We can now refine Theorem 1 in the following way.
\proclaim 
Theorem 3. The inequalities 
$$
\displaylines{
\quad 
{ (\omc_{n-1}(\Th) -1)\om_{n-1}(\Th) \over ((n-2)\omc_{n-1}(\Th) +1)\om_{n-1}(\Th) +(n-1)\omc_{n-1}(\Th)}
\le \om_0(\Th) \le 
\hfill\cr \hfill
{(1 - \omc_0(\Th))\om_{n-1}(\Th) -n +2 -\omc_0(\Th) \over n-1},\qquad
}
$$
hold for any point $\Th$ in $\bR^n$ with $1, \th_1, \dots , \th_n$ linearly independent over $\bQ$.

The above  estimate is stronger than (1), since 
$$
\omc_{n-1}(\Th) \ge n \and \omc_0(\Th) \ge {1\over n}
$$
by Dirichlet's Box Principle. 

Notice that Theorem 3 is optimal in dimension $n=2$. In that case, the uniform exponents $\omc_0(\Th)$
and $\omc_1(\Th)$ are linked by  Jarn\'\i k's equation [\JarD]
$$
\omc_0(\Th) = 1 -{1\over \omc_1(\Th)},  \leqno{(7)}
$$
so that  Theorem 3 reads equivalently
$$
{(\omc_1(\Th) -1)\om_1(\Th)\over \om_1(\Th) + \omc_1(\Th)}
\le \om_0(\Th) \le {\om_1(\Th) - \omc_1(\Th)+1 \over \omc_1(\Th)}. \leqno{(8)}
$$
It is established in [\LauB] that   the subset  of $(\bR_{>0}\cup \{+\infty\})^4$ made up  by all possible quadruples
$$
\Big(\om_1(\Th),\om_0(\Th),\omc_1(\Th),\omc_0(\Th)\Big), 
$$
when  $\Th$ ranges  over $\bR^2$ with $1,\th_1 , \th_2$ linearly independent over $\bQ$, is   essentially described  by the 
relations (7-8), together with the obvious lower bound $\omc_1(\Th)  \ge 2$. 
Thus, the transfer inequalities (8) cannot be sharpened for generic points $\Th$.

To end this introduction, let us ask for various extensions of the above results. First, Khintchine's Transference Principle 
has  naturally been extended to any system of real linear forms [\Dys ,  \SchB], and even more generally, to systems
of linear inequalities in an adelic context [\Bur].  We address the problem of splitting these further transfer relations through
adequate intermediate exponents, on the model of Theorem 2. That may  possibly  be achieved  by  improving the transfer 
results of [\SchA]. In an other direction,  one should also ask for some refined version   of Theorem 2, 
involving  uniform exponents as in Theorem 3. 

Transference  principles play  an important  role   in proving Schmidt's Subspace Theorem (see e.g. [\SchB]), which extends the classical 
Roth's theorem on the  rational approximations to an algebraic number [\Rot].  They also occur in 
  questions   of approximation to a real number $\xi$ by algebraic numbers of bounded
degree.  With regard to that topic, it might  be fruitful  to investigate more specifically transfer inequalities   between 
 exponents of the form  $\om_d(\xi, \dots , \xi^n)$,  
as well as their  uniform analogues $\omc_d(\xi, \dots, \xi^n)$  introduced (over the notation $\hat{w}_{d,1}(\Th)$ with $\Th =(\xi, \dots,\xi^n)$)
 in Definition 4 of [\BuLaC]. See the articles  [\DaSc, \LauA, \RoyA, \RoyC] and the reports [\BuLaC, \RoyB]
for  relating    those  exponents  with    various results on algebraic approximation.

\section
{2. Algebraic formulation}
We reformulate in this section Definition 2 in term of wedge product. One  relates  the exponent $\om_d(\Th)$ to integer solutions 
of a system of linear inequations, as in Definition 1.  This approach will enable us to employ standard arguments from the geometry of numbers.

First, we equip the real vector space $\bR^{n+1}$ with the usual dot product and extend it naturally to the Grassmann algebra 
$\Lambda(\bR^{n+1})$, by requiring that for any orthonormal basis $\{e_i\}_{0\le i\le n}$ of $\bR^{n+1}$, the family of  wedge products 
$$
e_{i_1}\wedge \dots \wedge e_{i_k} \, ;  \quad  0 \le i_1 < \dots < i_k \le n, 0 \le k\le n+1,
$$
is an orthonormal basis of $\Lambda(\bR^{n+1})$.
For any $\bX \in \Lambda(\bR^{n+1})$, we denote by $| \bX |$ its Euclidean norm.  Now, let $P$ and $Q$ be points in $\bP^n(\bR)$
with homogeneous coordinates ${\bx}$ and ${\by}$. The projective distance between $P$ and $Q$ is measured by the ratio
$$
{\rm d}(P,Q) ={ |{\bx} \wedge {\by} | \over | {\bx} | | {\by} |},
$$
which does not depend on the choices of ${\bf x}$ and ${\bf y}$.  Using Lagrange's identity
$$
| {\bx} \wedge {\by} |^2  + ( {\bx} \cdot {\by})^2 = | {\bx} |^2 | {\by}|^2,
$$
we see that ${\rm d}(P,Q)$  turns out to be  the sine of the acute angle determined in $\bR^{n+1}$ by the two lines $\bR{\bf x}$ and $\bR{\bf y}$.
See [\Phi , \RoyB] for further properties of the projective distance.

Now, let $L$ be a $d$-dimensional linear subvariety in $\bP^n(\bR)$. Write  $ L= \bP(V)$, where $V$ is the $(d+1)$-dimensional
 subspace of $\bR^{n+1}$ spanned by the homogeneous coordinates of the points of $L$.
  Select a basis $\{\bx_0, \dots , \bx_d\}$ of $V$ and put 
  $$
   \bX= \bx_0\wedge \dots \wedge \bx_d.
   $$
  Using the canonical basis $\{e_i\}_{0\le i\le n}$ of $\bR^{n+1}$, 
  we may identify $\Lambda^{d+1}(\bR^{n+1})$ with $\bR^{{n+1\choose d+1}}$.
  The multivector $\bX$ is then called a  system of {\it Pl\" ucker coordinates} of $L$ (or of $V$).  Note that $\bX$ is
   determined up to multiplication by a non-zero real number. It is known that the correspondence $L \mapsto \bX$ establishes a bijection between 
  the set of $d$-dimensional linear subvariety of $\bP^n(\bR)$ and the  set of non-zero decomposable multivectors
  in $\Lambda^{d+1}(\bR^{n+1})$, up to an  homothety.

\proclaim
Lemma 1. Let $\Th$ be a point in $\bP^n(\bR)$ with homogeneous coordinates ${\by}$ and let $L$ be a linear subvariety of $\bP^n(\bR)$ with
Pl\" ucker coordinates $\bX$. Then
$$
{\rm d}(\Th , L) = { | {\by} \wedge \bX| \over | {\by} | | \bX| }.
$$

\proof
Write $L =\bP(V)$ as above.  If $\by$ is orthogonal to $V$, we have
$$
{\rm d}(\Th , L) = 1 = { | {\by} \wedge \bX| \over | {\by} | | \bX|} .
$$
Otherwise,  denote by $\by'$ the orthogonal projection on $V$ of the vector $\by$.  The minimal angle 
$\langle \bR\by , \bR\bx \rangle$, when $\bx$ ranges along $V\setminus \{ 0\}$, is  clearly reached for
$\bx = \by'$. Therefore
$$
{\rm d} (\Th , L) = { | \by  \wedge \by' | \over | \by | | \by' |}.
$$
We may assume without loss of generality that $\by\notin V$.  Let $\{e_i\}_{1\le i \le d+1}$  be an orthonormal basis of $V$, which
we complete by an orthogonal unitary vector $e_0$ to get an orthonormal basis $\{e_i\}_{0\le i \le d+1}$ of $\bR \by\oplus V$.
Choose $\bX = e_1\wedge \dots \wedge e_{d+1}$ as a system of Pl\" ucker coordinates of $V$.
Then $| \bX |=1$.  Write now  $ \by = a e_0 + \by'$ for some  $a\in \bR$. We have 
$$
|  \by\wedge \by' | = | a| | e_0  \wedge \by' | = | a| | \by' | \and 
| \by \wedge \bX | = | a| | e_0 \wedge e_1 \wedge \dots \wedge e_{d+1} | = | a |,
$$
so that
$$
{\rm d}(\Th , L) ={ | \by \wedge \by' | \over | \by | | \by' |}= { | a| \over | \by |} = { | \by \wedge \bX  | \over | \by | | \bX |},
$$
as required. \cqfd

\bigskip

Suppose now that $L$ is rational over $\bQ$, or equivalently that $V$ can be generated by points $\bx_0, \dots , \bx_d$
belonging to $\bQ^{n+1}$.  Then $\bX= \bx_0 \wedge \dots \wedge \bx_d$ has rational coordinates in $\bR^{{n+1\choose d+1}}$.
We define  $H(L)$ as the Weil height of the ${n+1\choose d+1}$-tuple $\bX$.
Let us  indicate  an useful interpretation of the height $H(L)$  in term of determinant of a lattice. 

\proclaim
Lemma 2. The group $V \cap \bZ^{n+1}$ is a lattice in
$V$. 
For any  $\bZ$-basis $\{\bx_0, \dots , \bx_d\}$  of $V \cap \bZ^{n+1}$, put $\bX= \bx_0\wedge\dots \wedge \bx_d$. 
Then
$$
\det (V\cap \bZ^{n+1}) = | \bX | = H(L).
$$

The notation $\det(V\cap \bZ^{n+1})$ indicates here the volume, with respect to the induced  Euclidean norm on $V$,
 of the unit parallelepiped constructed  on any such $\bZ$-basis $\{\bx_0, \dots , \bx_d\}$.
 
 \proof 
 The $\bZ$-module $V \cap \bZ^{n+1}$ has rank $d+1$,  since $V$ is rational over $\bQ$.
 The multivector $\bX$ is clearly a system of Pl\" ucker coordinates of $L$, and  $\bX$ is primitive in 
 the group $\Lambda^{d+1}(\bZ^{n+1})$ since $\{\bx_0, \dots , \bx_d\}$ is a $\bZ$-basis of $V \cap \bZ^{n+1}$.  
 Therefore $| \bX | = H(L)$. Moreover
 $$
 \det (V\cap \bZ^{n+1})= \sqrt{\det \Big( \bx_i\cdot\bx_j\Big)_{ {0\le i \le d\atop 0  \le j\le d}}} = | \bX |,
 $$
 by Cauchy-Binet formula.
 We refer to  Theorem 1 of [\SchA] for more details  and for an extension of Lemma 2 to number fields.
 \cqfd

\bigskip
Lemmas 1 and 2 enable us to handle  more easily  the quantities  $ \om_d(\Th)$ through the following alternative

\proclaim
 Definition 4. Let $\by$ be homogeneous coordinates of the point $\Th$.
 For any $d$ with $0 \le d \le n-1$, the exponent $\om_d(\Th)$ is the supremum of the real numbers $\om$
for which there exist infinitely many integer decomposable multivectors $\bX \in \Lambda^{d+1}(\bZ^{n+1})$ such that
$$
| \by \wedge \bX | \le | \bX |^{-\om}.
$$

 For any $d$-dimensional rational linear subvariety $L$, select
a system of Pl\" ucker coordinates $\bX$ as in Lemma 2. Then
$$
H(L) = | \bX | \and {\rm d} (\Th , L) = | \by |^{-1} | \bX |^{-1} | \by \wedge \bX |.
$$
Now, the equivalence of Definitions 2 and 4  is clear.

\bigskip
\noindent
{\bf Remark. } We should have obtained the same exponent $\om_d(\Th)$,  when dropping in Definition 4 the assumption that $\bX$
is decomposable in $\Lambda^{d+1}(\bZ^{n+1})$. We shall not use this property, which can be deduced from Mahler's theory of
compound convex bodies [\Mah, \SchB]. To that purpose, observe that for fixed positive real numbers $Y  <X$, 
the convex body defined as the set of $\bX$ in $\Lambda^{d+1}(\bR^{n+1})$ satisfying the linear inequations
$$
| \bX | \le Y^d X  \and | \by \wedge \bX | \le Y^{d+1},
$$
compares with the $(d+1)$-th compound of the convex body defined in $\bR^{n+1}$ by 
$$
| \bx | \le X \and | \by \wedge \bx | \le Y.
$$
Further applications to transfer inequalities will be given elsewhere.

\bigskip

\section
{3. Proof of the Going-up inequality.}
We are now able to prove inequality (2). Fix homogeneous coordinates $\by$ of $\Th$ and follow  Definition 4.
Let $\bX$ be a decomposable multivector in $\Lambda^{d+1}(\bZ^{n+1})$ such that
$$
| \by \wedge \bX | \le | \bX |^{-\om}.
$$
We may suppose without loss of generality that $\bX$ is primitive in the group $\Lambda^{d+1}(\bZ^{n+1})$. Then $\bX$ is a system of 
Pl\" ucker coordinates of some rational $(d+1)$-dimensional subspace $V \subset \bR^{n+1}$.  By Lemma 2, we may write
$\bX = \bx_0 \wedge \dots \wedge \bx_d$ for some $\bZ$-basis $\{\bx_0, \dots , \bx_d\}$ of the lattice $V\cap \bZ^{n+1}$. 
Specifically,  we know from Lemma 2 that 
$$
\det ( V \cap \bZ^{n+1}) = | \bX |.
$$
Let $W$ be the orthogonal complement  to  $V$ in $\bR^{n+1}$. Therefore, $W$ is a rational subspace of $\bR^{n+1}$ with  dimension $n-d$. 
Denote by $\Lambda$ the orthogonal projection on $W$ of the lattice $\bZ^{n+1}$. Then, $\Lambda$ is a lattice in $W$ with determinant
$$
\det (\Lambda ) = { \det ( \bZ^{n+1} ) \over \det ( V \cap \bZ^{n+1})} = { 1 \over | \bX |}.
$$
On the other hand, the Euclidean ball in $W$
$$
\{ \bz \in W \, ; \quad | \bz | \le R\},
$$ 
centered at the origin of $W$ with radius $R$, has volume $v_{n-d}R^{n-d}$, where $v_k = \pi^k/\Gamma(1 +k/2)$ denotes the
volume of the unit Euclidean ball in $\bR^k$. Choosing now
$$
R = 2 v_{n-d}^{-1/(n-d)} | \bX |^{-1/(n-d)},
$$
Minkowski's convex body theorem shows that there exists a non-zero element $\xi \in \Lambda$ with norm
$$
| \xi | \le R.
$$
Lift up now $\xi$ into an element $\bx \in \bZ^{n+1}$ whose  orthogonal projection  on $W$ is $\xi$, and put
$$
\bX' = \bx \wedge \bX = \xi \wedge \bX.
$$
Then $\bX'$ is a decomposable multivector in $\Lambda^{d+2}(\bZ^{n+1})$.  Making use of  Hadamard's inequality, 
we first bound from above its norm
$$
| \bX' | = | \xi | | \bX| \le R | \bX | \ll | \bX | ^{(n-d-1)/(n-d)}.
$$ 
A  second  use of Hadamard's inequality enables us to bound the wedge product
$$
| \by \wedge \bX' | = | \by \wedge \xi \wedge \bX | \le  |  \xi | | \by \wedge \bX |
\le   R | \bX |^{-\om} \ll | \bX |^{-(\om + 1/(n-d))}.
$$
Combining the two last inequalities, we obtain
$$
| \by \wedge \bX' | \ll | \bX' |^{ -( (n-d)\om +1)/(n-d-1)}.
$$
Taking now $\om$ arbitrarily close to $\om_d(\Th)$, we have established the lower bound (2).

\vfill\eject

\centerline{\bf References }

\vskip 1cm

\item{[\BuLaB]}
Y. Bugeaud and M. Laurent,
{\it On exponents of homogeneous and inhomogeneous Diophantine Approximation}, 
 Moscow Mathematical Journal, 5, 4 (2005), 747--766.

\item{[\BuLaC]}
Y. Bugeaud and M. Laurent,
{\it Exponents of Diophantine Approximation}, 
to appear in the Publ. of the Ennio de Giorgi Math. Institute, Pisa.

\item{[\Bur]}
E.  Burger,
{\it Homogeneous Diophantine approximations in $S$-integers}, 
Pacific J. Math.  152 (1992), 211--253.

\item{[\Cas]}
        J. W. S. Cassels,
An introduction to Diophantine Approximation,
Cambridge Tracts in Math. and Math. Phys., vol. 99, Cambridge
University Press, 1957.

\item{[\DaSc]}
        {H. Davenport and W. M. Schmidt},
{\it Approximation to real numbers by
algebraic integers}, Acta Arith. {15} (1969), 393--416.

\item{[\Dys]}
F. J. Dyson,
{\it On simultaneous Diophantine approximations},
Proc. London Math. Soc. 49 (1947), 409--420.

\item{[\Gro ]} 
A. Groshev, {\it  A theorem on a system of linear forms}, Dokl. Akad Nauk SSR 19 (1938), 151-152 (in Russian).

\item{[\JarA]}
V. Jarn\'\i k,
{\it \"Uber ein Satz von A. Khintchine},
Pr\'ace Mat.-Fiz. 43 (1935), 1--16.

\item{[\JarC]}
V. Jarn\'\i k,
{\it \"Uber ein Satz von A. Khintchine, 2. Mitteilung},
Acta Arith.  2  (1936), 1--22.

\item{[\JarD]}
V. Jarn\'\i k,
{\it Zum Khintchineschen ``\"Ubertragungssatz''},
Trav. Inst. Math. Tbilissi 3 (1938), 193--212.

\item{[\KhiA]}
A. Ya. Khintchine,
{\it Einige S\" atze \" uber Kettenbr\" uche, mit Anwendungen auf die Theorie der Diophantische Aproximationen},
Math. Ann. 92 (1924), 115--125.

\item{[\KhiB]}
A. Ya. Khintchine,
{\it \"Uber eine Klasse linearer diophantischer Approximationen},
Rendiconti Circ. Mat. Palermo 50 (1926), 170--195.

\item{[\LauA]}
        M. Laurent,
{\it Simultaneous rational approximation to the successive powers
of a real number},
Indag. Math. 11 (2003), 45--53.

\item{[\LauB]}
M. Laurent,
{\it Exponents of Diophantine Approximation in dimension two},
to appear in the Canad. Journal of Math..

\item{[\Mah]}
          K. Mahler,
{\it On compound convex bodies, I},
Proc. London Math. Soc. (3)  {5} (1955), 358--379.

\item{[\Phi]}
          P. Philippon,
{\it Crit\`eres pour l'ind\'ependance alg\'ebrique},
Inst. Hautes \'Etudes Sci.  Publ. Math.    {64} (1986), 5--52.

\item{[\Rot]}
          K. F.  Roth,
{\it Rational approximations to algebraic numbers},
Mathematika 2 (1955), 1--20.

\item{[\RoyA]}
        D. Roy,
{\it Approximation to real numbers by cubic algebraic numbers, I},
Proc. London Math. Soc. 88 (2004), 42--62.

\item{[\RoyC]}
        D. Roy,
{\it Approximation to real numbers by cubic algebraic numbers, II},
Annals of Math. 158 (2003), 1081--1087.

\item{[\RoyB]}
          D. Roy,
{\it Diophantine approximation in small degree},
Number Theory, C.R.M. Proc. Lecture Notes 36 (2004), 269--285.

\item{[\SchA]}
W. M. Schmidt,
{\it On heights of algebraic subspaces and Diophantine approximations},
Ann. Math.  85 (1967), 430--472.

\item{[\SchB]}
W. M. Schmidt, Diophantine Approximation,  Lecture Notes in Mathematics, vol. 785, 
Springer, Berlin, 1980.

\bigskip

\noindent   {Michel LAURENT}

\noindent 
{Institut de Math\'ematiques de Luminy}

\noindent 
{C.N.R.S. -  U.M.R. 6206 - case 907}

\noindent       {163, avenue de Luminy}

\noindent 
{13288 MARSEILLE CEDEX 9  (FRANCE)}

\noindent 
{\hbox{\tt laurent@iml.univ-mrs.fr}}

\bye